\documentclass[reqno, 11pt]{amsart}
\usepackage[utf8]{inputenc}
\usepackage[T1]{fontenc}
\usepackage{lmodern}
\usepackage{epsfig,amsmath, amsthm ,amsfonts,latexsym,stmaryrd, amssymb}
\usepackage[abs]{overpic}
\usepackage[dvipsnames]{xcolor}
\usepackage[linktocpage=true]{hyperref}
\usepackage{caption}
\usepackage{subcaption}
\usepackage{dsfont}
\usepackage{mathtools}
\usepackage{tikz-cd}
\usetikzlibrary{shapes,arrows}
\usepackage{comment}
\usepackage{enumitem}
\usepackage{marginnote}
\usepackage[all,cmtip]{xy}
\usepackage[backend=biber, maxnames=4, style=alphabetic]{biblatex}
\hypersetup{
  colorlinks,
  citecolor=teal,
  linkcolor=purple,
  urlcolor=teal}
\theoremstyle{plain}

\newtheorem{dummy}{anything}[section]

\newtheorem{lemma}[dummy]{Lemma}

\newtheorem{corollary}[dummy]{Corollary}

\newtheorem{proposition}[dummy]{Proposition}

\theoremstyle{definition}
\newtheorem{definition}[dummy]{Definition}

\theoremstyle{remark}
\DeclareFontFamily{U}{mathx}{}
\DeclareFontShape{U}{mathx}{m}{n}{<-> mathx10}{}
\DeclareSymbolFont{mathx}{U}{mathx}{m}{n}
\DeclareMathAccent{\widecheck}{0}{mathx}{"71}
\begin{document}

\title[On cone structures]{Cone structures from a dynamical and probabilistic viewpoint}

\author{Agustin Moreno}

\address[A.\ Moreno]{Institut f\"ur Mathematik \\ Universit\"at Heidelberg \\ Germany}

\email{\href{mailto:agustin.moreno2191@gmail.com}{agustin.moreno2191@gmail.com}}

\date{\today}

\begin{abstract}
    The goal of this note is to explore, from a geometric and probabilistic point of view, the dynamics of cone structures adapted to open book decompositions. This is inspired by the picture which arises in the study of the circular restricted three body problem (CR3BP). This yields geometric obstructions to reaching a point from another point in the phase space of the CR3BP.
\end{abstract}

\maketitle

\section{Introduction}

The goal of this note is to explore, from a geometric and probabilistic point of view, the dynamics of cone structures (i.e.\ distributions of conical sets) adapted to open book decompositions. This is inspired by the picture which arises in the study of the circular restricted three body problem, see \cite{M23}. A problem that inspired this investigation is that of understanding which points in configuration space can be reached from a predetermined point.

\medskip

\textbf{Acknowledgements.} The author would like to express gratitude to Mischa Gromov, Dennis Sullivan, Helmut Hofer, Urs Frauenfelder, Alberto Abbondandolo, Peter Albers, Otto van Koert, Rohil Prasad, Connor Jackman.

\medskip

\textbf{Setup.} A \emph{concrete} open book decomposition $\theta: M \backslash B \rightarrow S^1$ on a contact odd-dimensional manifold $(M,\xi)$ is, by definition, a fibration which coincides with the angle coordinate on a choice of collar neighborhood $B\times \mathbb{D}^2$ for a codimension-2 closed submanifold $B\subset M$ (the \emph{binding}). We assume that it supports $\xi$ in the sense of Giroux. This means that there is a contact form $\alpha$ for $\xi$, a \emph{Giroux form}, such that $\alpha\vert_B$ is contact, and $d\alpha$ is positively symplectic on the fibers of $\theta$; equivalently, the Reeb flow of $\alpha$ has $B$ as an invariant subset, and it is positively transverse to each fiber. We denote the \emph{$\varphi$-page} by $P_\varphi=\overline{\theta^{-1}}(\varphi)$ for $\varphi \in S^1$, and we also use the abstract notation $M=\mathbf{OB}(P,\phi)$, where $P$ is the abstract page (the closure of the typical fiber of $\theta$) with $\partial P=B$, and $\phi$ is the symplectic monodromy. 

In \cite{M23}, the author considered the dynamics of the circular, restricted three-body problem (or CR3BP for short). From \cite{MvK}, the (low energy, near primary) dynamics of this problem can be viewed as a $5$-dimensional situation as above. Whenever the planar dynamics is convex/dynamically convex, as follows from \cite{HSW}, combined with Theorem 1 in \cite{MvK}, one can construct a ``shadow'' dynamics on the three-dimensional sphere, see \cite{M23}. This consists of a Reeb dynamics on a moduli space of pseudo-holomorphic curves, which is obtained by a suitable integration procedure, using a finite energy foliation of the original manifold. This induced dynamics is given by an averaged contact form on the moduli space.
In what follows, we will then focus on the case $M=S^3=\mathbf{OB}(D^2,\mathds{1})$, with the trivial open book whose pages are disks (and therefore the monodromy is the identity), while the contact structure is the standard one. 

Moreover, by projecting the original Reeb field on the $5$-manifold to $S^3$, the author constructed a cone structure, i.e.\ a choice of conical set $C=\{C_p\}_{p \in S^3}$ depending smoothly on $p$. Here, a conical set is a subset of the tangent bundle to $M$, invariant under the natural $\mathbb R^+$ scaling action. The cone structure is trivial at a point $p$ if $C_p=\{0\}$. The cone structure $C$ which arises on $S^3$ is $2$-dimensional, i.e.\ it is the image of a map with a two-dimensional domain. The orbits of the original flow project, by construction, to smooth curves tangent to the cone structure. With the orientation given, the tangent plane $TP$ to the page determines a half-space of positive directions. A \emph{base} for a cone inside a tangent space is a section of the cone, i.e.\ a continuous choice of vector in each direction tangent to the cone; a base for a cone structure is a continously varying choice of bases for each tangent space, see e.g.\ \cite{S76}. A \emph{based} cone structure is a cone structure for which a continuously varying base has been chosen. Such a based cone structure naturally arises in the picture for the circular, restricted three-body problem. There is a natural Riemannian metric on $S^3$ associated to a suitable almost complex structure and the averaged contact form, which will be implicitly used in what follows to measure angles and distances. In this context, we have the following.

\begin{figure}
    \centering
    \includegraphics[width=0.5 \linewidth]{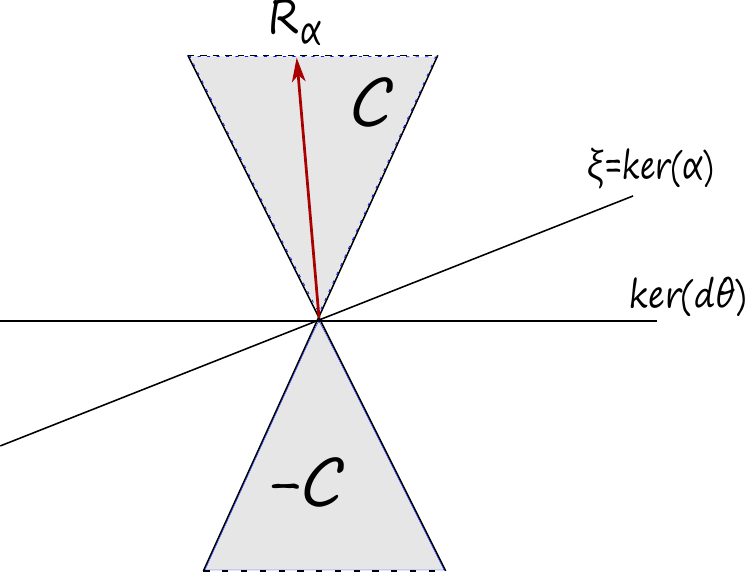}
    \caption{A cone structure adapted to an open book.}
    \label{fig:adapted}
\end{figure}

\begin{definition}\label{def:stronglyadapted}
Consider an everywhere non-trivial cone structure $C$ on a manifold $M$, where $M$ is endowed with an open book $\theta: M\backslash B\rightarrow S^1$. We say that $C$ is adapted to $\theta$ if
\begin{enumerate}
    \item[(1)] $C\vert_B\subset TB$;
    \item[(2)] $d\theta$ is a section of $C\vert_{M\backslash B}$; 
\end{enumerate}
and if there exists a Giroux form $\alpha$ for the open book such that
\begin{enumerate}
    \item [(3)] $\alpha$ is a section for $C$;
    \item [(4)] The Reeb vector field $R_\alpha$ is interior to $C$.
\end{enumerate}
Here, a section for $C$ is a $1$-form which is strictly positive on non-zero vectors of $C$. See Figure \ref{fig:adapted}. Moreover, a vector is interior to a conical set if the smallest disk which covers some base of the cone contains the vector.
\end{definition}

\section{Dynamics of cone structures}

Consider a cone structure $C$ adapted to an open book, on a contact manifold $(M,\xi)$. We will assume, for simplicity, that $M=S^3$. We say that a $C^1$ curve $\alpha$ is a trajectory for the cone structure, if it is tangent to the cone at every point. Given a point $p$, its future time $t$ region consists of all those points which are reachable in time $t$ by curves tangent to the cone, parametrized to have unit speed. We have the following lemma, whose proof is a simple exercise.

\begin{figure}
    \centering
    \includegraphics[width=0.37\linewidth]{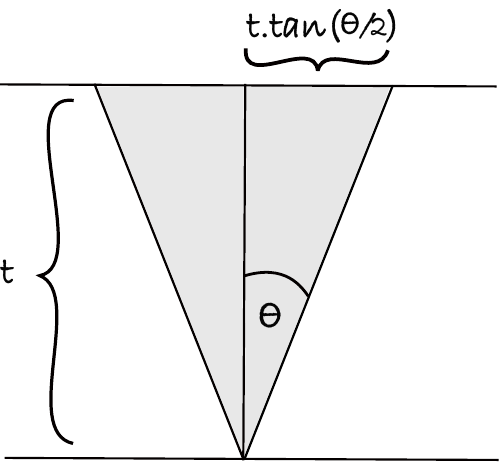}
    \caption{Orbits emanating from a point can reach, in time $t$, points lying in a disk of radius $t\cdot tan(\theta/2)$, where $\theta$ is the inner angle of the cone.}
    \label{fig:cono}
\end{figure}

\begin{lemma}\label{lem:angle}
    In the upper half-plane $\mathbb R^3_+$, consider a cone structure with constant inner angle $\theta$.\footnote{Here, the angle is measured between an edge of the cone, and a vertical line, see Figure \ref{fig:cono}.} Then the future time-$t$ region is a disk of radius $t\cdot \mbox{tan}(\theta/2)$. $\hfill$ $\square$
\end{lemma}

Heuristically, the above lemma gives a geometric obstruction to connecting the starting point to another arbitrary point. Explicitly, the starting point cannot be connected to a given putative end point, if the latter lies outside a disk of radius as given by the lemma. As a sanity check, note that $\theta=0$ corresponds to zero radius, whereas $\theta = \pm \pi$ corresponds to infinite radius, and $\theta = \pi/4$ corresponds to radius $1$.

Fix a page $P=P_0$ of the open book. Consider a disk $A=A(p,r)$ lying inside $P$, of radius $r$, centered at a point $p$. We assume that the cone structure $C$ on the sphere has constant inner angle $\theta$ (we can always take the smallest algebraic cone\footnote{Here, by \emph{algebraic}, we mean a cone which admits a disk of radius $r$ as a base, for some $r>0$.} with constant inner angle containing $C$). We also assume that the disk-like page of the open book has radius $1$, for simplicity. Consider the corresponding $t$-page $P_t$, and a measurable subset $B\subset P_t$. We wish to estimate the probability $p^t_{A,B}$ that points from $A$ arrive to $B$ (in time $t$).\footnote{The time is measured with respect to the open book, i.e.\ arriving at the $t$-page $P_t$ from the $0$-page $P_0$ is considered time $t$, independent on the orbit. In other words, the ``time-keeper'' is the Hopf fibration.} Let $\mu_t$ be the Lebesgue measure on $P_t$ induced by the averaged contact form.

\begin{lemma}\label{eq:prob}
We have
$$
p^t_{A,B}= \mu(B\cap D_{A,\theta}^t),
$$
where $D_{A,\theta}^t$ is a disk of radius $t \cdot \mbox{tan}(\theta/2) \mu(A)$, centered at the image $p_t$ of $p$ under the time-$t$ Hopf flow, with $\mu(A)=\pi r^2$.
\end{lemma}

Note that if $\theta \geq \pi/4$, then there is no geometric obstruction to reaching the set $B$. This lemma, whose proof follows by integration of the result in Lemma \ref{lem:angle}, can be adapted to the case where $A$ is an arbitrary measurable set. 

\medskip

\textbf{Cone structures as measures of integrability.} The degenerate case when the cone structure $C$ gives a ray at each point is the case of a (reparametrization class of a) flow adapted to the open book. Heuristically, this case can be thought of as the ``integrable'' case, e.g.\ the Hopf flow, which within the context of the CR3BP, is the shadow dynamics of the integrable limit cases (i.e.\ the Kepler problem and its rotating version, see e.g.\ \cite{MvK22}). With this in mind, and denoting by $\varphi: M \rightarrow \mathbb R$ the function which at each point $p$ gives the inner angle of the smallest algebraic cone containing $C_p$, one can define two invariants, the \emph{mean measure} of integrability
$$
\mathcal{I}_m(C)=\int_M \varphi \cdot d\mbox{vol},
$$
and the \emph{max measure} of integrability
$$
\mathcal{I}_M(C) = \max_M \varphi.
$$
By construction, both measures (invariants of the cone structure) vanish for the case of a flow, e.g.\ for the Hopf flow. 

Note that, by construction, the invariant $\mathcal{I}_M(C)$ is the inner angle of a cone structure with constant inner angle containing $C$. Therefore the result of Lemma (\ref{eq:prob}) applies, and yields an upper bound, as follows.

\begin{corollary} Take a disk $A\subset P=P_0$ centered at $p\in P,$ of radius $r>0$, and a measurable subset $B\subset P_t$ for $t>0$. Then the probability of going from $A$ to $B$ in time $t$ via trajectories tangent to the cone admits a bound:
$$
p^t_{A,B}\leq \mu(   B   \cap D^t_{A,\mathcal{I}_M(C)})
$$
where $D^t_{A,\mathcal{I}_M(C)}$ is as in Lemma \ref{eq:prob}, with $\theta = \mathcal{I}_M(C)$.
\end{corollary}

\begin{figure}
    \centering
    \includegraphics[width=0.5\linewidth]{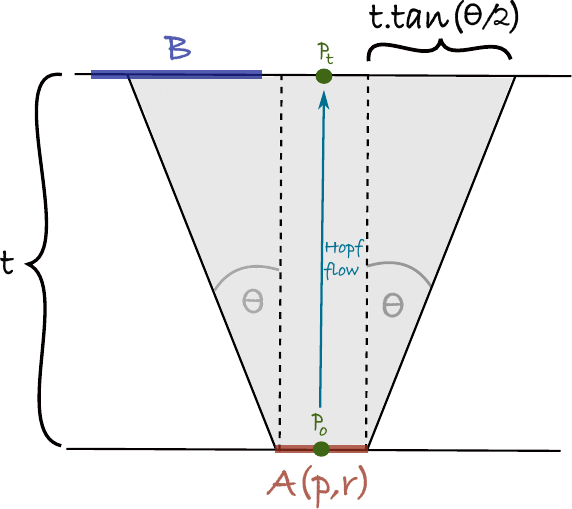}
    \caption{The geometric measure theory of constant cone structures.}
    \label{fig:placeholder}
\end{figure}

\section{Invariants for cone structures}

Let $\mathcal{P}_C$ be the set of $C^1$-trajectories tangent and/or interior to the cone structure, which start at $P_0$, and end at $P_0$, and which come with their own time parametrization. In other words, there is a map
$$
\tau: \mathcal{P}_C \rightarrow \mathbb R,
$$
for which $\tau(\gamma)$ is the first return time to $P_0$ of the curve $\gamma$. A \emph{section} of $\mathcal{P}_C$ is a choice $\Gamma=\{\gamma_p\}_{p \in P_0\backslash B}$ of curves in $\mathcal{P}_C$, such that the curve $\gamma_p$ starts at $p$, and it depends continuously on $p$ (lying away from the binding $B$, where the cone structure is assumed to degenerate to a line). 

Given a measurable set $A$, and a section $\Gamma$, we denote by 
$$\tau(\Gamma): S^3 \rightarrow \mathbb R$$ the (measurable) map that associates to $p \in S^3$, the return time of the unique curve in $\Gamma$ passing through $p$. We can then define the associated \emph{Calabi invariant}, given by
$$
\mbox{CAL}_\Gamma(A)=\int_A \tau(\Gamma)d\mu.
$$
In other words, this is the expected value of the return time, for orbits of $\Gamma$ passing through $A$. This assumes that the Poincar\'e return map extend continuously to the boundary, in order to have a well-defined, integrable return time function $\tau(\Gamma)$. From Appendix C in \cite{MvK}, for the case where $\Gamma$ consists of the Reeb flow orbits for the shadow dynamics, we gather that
$$
\mbox{CAL}_\Gamma(P)=\mbox{vol}(M),
$$
is the contact volume of the ambient manifold. In general, $\mbox{CAL}_\Gamma(A)$ agrees with the volume of the suspension of $A$.

If $\tau_n(\Gamma) :\mathcal{P}_C \rightarrow \mathbb R$ represents the $n$-th return time function (i.e.\ the time it takes for the orbit to return $n$ times to $P$), with respect to an arbitrary section $\Gamma$, it would be interesting to provide information on the growth rate of the limit
$$
\lim_{n\rightarrow +\infty} \mbox{CAL}^n_\Gamma(A) = \lim_{n\rightarrow +\infty} \int_{A_n} \tau_n(\Gamma)d\mu,
$$
where $A_n = \tau_n(A)$. Note that this limit equals $\mbox{vol}(M)$ in the above case, where $\Gamma$ consists of orbits of a Reeb flow.

%From now on, we assume for simplicity that we have normalized so that the total mass $\int_P \tau(\Gamma)d\mu$ is equal to $1$. 

\medskip

\textbf{Connection to QFT.} It is well-known that, in the Quantum Theory of Fields, the Feynmann path integral is ill-defined. In our setting, there is no well-defined measure on $\mathcal{P}_C$. However, given a section of $\mathcal{P}_C$, this allows to define a measure on the section, by identifying it with the interior of one of the pages of the open book. We can then consider a uniform distribution which samples points from the page $P=P_0$ (which we identify with a disk in the plane), which then yields the same distribution on the section. Note that its mean, and its variance are: 

$$
E(\Gamma)=\int_P\frac{x}{\mu(P)}d\mu(x) = \int_P x d\mu(x) ,
$$
is the center of mass of the page $P$ (i.e.\ the ``origin'' of the disk), and 
$$
\mbox{Var}(\Gamma)=\frac{1}{\mu(P)}\int\left(x-\int_P x d\mu(x)\right)^2d\mu(x)= \frac{\mu(P)^2}{12}=\frac{1}{12},$$
where we used that $\mu(P)=1$.

\section{Brownian motion for cone structures}

We will now introduce randomness to the setup. Consider a $3$-dimensional Wiener process $W_t=(X_t,Y_t,Z_t)\in \mathbb R^3$, and an associated three-dimensional stochastic process $V_t=(X'_t,Y'_t,Z'_t)$, given by the stochastic differential equation
\begin{equation}\label{eq:E_Gamma}
    dV_t = \sigma dW_t + \mu \tau(\Gamma) dt.
    \tag{$E_\Gamma$}
\end{equation}
Here, $\sigma: S^3 \rightarrow \mathbb R$ is the \emph{volatility}, $\mu>0$ is a constant drift term of the form $(0,0,\mu_3)$, and $\tau(\Gamma)$ is the measurable return time function encountered above, associated to a section $\Gamma$. Intuitively speaking, the process $V_t$ has a tendency to move forwards along the $z$-axis, while other components have no tendency to rise or fall. 

This type of process can also be adapted to an open book with an adapted cone structure (e.g.\ for our running example on $S^3$), by imposing that the Brownian motion is interior to the cone structure, and has a positive drift in the open book angle direction. In this situation, given a section $\Gamma$ for $\mathcal{P}_C$, the function $\tau(\Gamma)$ can be interpreted as the average return time to a fixed page, of a curve in the section. From standard results of Brownian motion, we obtain the following.

\begin{proposition}
Given a point $p\in P_0$, an open subset $U\subset P_0$, and a section $\Gamma$ of $\mathcal{P}_C$, there is probability $1$ that there is a development $V_t \in S^3$ of the above random process (\ref{eq:E_Gamma}) which starts at $p$, and ends at $U$, in some integer time $n>0$. The average time for such $n$ is, however, infinite.
\end{proposition}

Note that in a deterministic situation (e.g.\ coming from a Hamiltonian system), the above proposition is compatible with Poincar\'e recurrence. We will then call the above, \emph{probabilistic Poincar\'e recurrence}. Heuristically, from the perspective of quantum mechanics, and in particular, the Feynmann path integral, one can interpret the smooth curves in $\mathcal{P}_C$ as the ``classical solutions'', and their probabilistic version (the instantiations of the process $V_t$) as the ``quantum oscillations around the classical solution''.

\printbibliography

\end{document}